\tolerance=20000
\overfullrule=0pt
\magnification=\magstep1
\mathsurround=2pt
\footline={\ifnum\pageno=1{\hfil}\else{\hss\tenrm\folio\hss}\fi}
\font\bigbf=cmbx10 scaled \magstep1

\font\biggtenex=cmex10 scaled \magstep2
\font\tenmsy=msym10
\font\sevenmsy=msym7
\font\fivemsy=msym5
\newfam\msyfam
\textfont\msyfam=\tenmsy
\scriptfont\msyfam=\sevenmsy
\scriptscriptfont\msyfam=\fivemsy
\def\Bbb#1{\fam\msyfam#1}
\baselineskip=16pt
\lineskip=3pt
\lineskiplimit=3pt
\def\cW{{\cal W}}
\def\R{{\Bbb R}}
\def\topsim#1{\ \smash{\mathop{\sim}\limits^{#1}}\ }
\def\ttopsim#1{\ {\buildrel{#1}\over\sim}\ }
\def\bigk{\mathop{\biggtenex K}}

\def\square{\vbox{
      \hrule height 0.4pt
      \hbox{\vrule width 0.4pt height 5.5pt \kern 5.5pt \vrule width 0.4pt}
      \hrule height 0.4pt}}

\def\Pn1{P_{n}^{(1)} (z)}
\def\Pn0{P_{n}^{(0)} (z)}

\def\tphin{\> {}_{10}\phi_9}
\def\ephis{\> {}_8\phi_7}

\def\fphit{\> {}_4\phi_3}

\def\bibline{\ {\vrule width.5in height .03pt depth .05pt}\ }
\rm
\vglue .5in
\centerline {\bigbf The Last of the Hypergeometric Continued Fractions
 \footnote{*}{{\rm Research partially
supported by NSERC (Canada).A talk prepared for the Internationl Conference
on Mathematical Analysis and Signal Processing, Cairo, Egypt, Jan.
3-9, 1994.}}}
\vglue .25in
\centerline {David R. Masson}
\bigskip
\baselineskip=12pt
\centerline {Department of Mathematics}
\centerline {University of Toronto}
\centerline {Toronto Ontario M5S 1A1}
\centerline {Canada}
\bigskip
\vglue .75in

\baselineskip=14pt
\midinsert\narrower\narrower
\centerline {\bf Abstract}
\medskip

A contiguous relation for complementry pairs of very well poised
balanced ${}_{10}\phi_9$ basic hypergeometric functions is used to derive
an explict expression for the associated continued fraction.
This generalizes the continued fraction results associated with both  Ramanujan's
Entry 40 and Askey-Wilson polynomials which can be recovered as limits.
Associated with our continued fraction results there are 
systems of 
biorthogonal rational functions that have yet to be derived.
\bigskip

\medskip

\endinsert
\vfill\eject

\baselineskip=18pt
\beginsection  1. Introduction.

In this paper we examine the level above the $ q $-Askey-Wilson 
polynomial case [1] using the methods of [4]. All of our parameters
are, in general, complex. We essentially follow the notation in Gasper and
Rahman [2], except that we omit the designation $ q $ for the base
in the q-shifted factorials and basic hypergeometric functions.
Thus we have, with $ | q | < 1 $,
$$
(a)_0:=1, \quad (a)_n:=\prod_{j=1}^n(1-aq^{j-1}),\quad n>0, \hbox
{ or }  n=\infty,
$$
$$
(a_1, a_2, \cdots, a_k)_n:= \prod_{j=1}^k (a_j)_n .
$$
The basic hypergeometric series $ _{r+1} \phi_r $ is
$$
_{r+1} \phi_r\left( \matrix{ a_1,&a_2,&\cdots,&a_{r+1} \cr
b_1,& b_2,&\cdots
,&b_r \cr}
\quad;z\right):=\sum_{n=0}^{\infty} {(a_1,a_2, \cdots,a
_{r+1})_n \over(b_1,b_2,\cdots,b_r,q)_n}\,z^n\quad.
$$
We denote a very well poised, balanced ${}_{10}\phi_9 $ as
$$ \phi \ = \ \phi(a; b,c,d,e,f,g,h;q)
$$
$$ \ := \  
\tphin \left({ {a, q\sqrt a , -q \sqrt a , b, c, d, e, f, g, h  }
\atop {\displaystyle\sqrt a, -\sqrt a, {aq\over b}, {aq\over c}, \ldots ,
{aq\over h}}} ; \ q\right) , bcdefgh=a^3q^2 .
$$
A complementary pair of very well poised, balanced ${}_{10}\phi_9$'s is
$$
\Phi^{(b)}(a;b,c,d,e,f,g,h;q) \ := \  \phi(a;b,c,d,e,f,g,h;q)
$$
$$+{(aq,{b\over a},c,d,e,f,g,h,{bq\over c},{bq\over d},{bq\over e},
{bq\over f},{bq\over g},{bq\over h})_\infty \over ({b^2q\over a},{a\over b},
{aq\over c},{aq\over d},{aq\over e},{aq\over f},{aq\over g},{aq\over h},
{bc\over a},{bd\over a},{be\over a},{bf\over a},{bg\over a},{bh\over a}
)_\infty}
$$
$$ \times \phi(b^2/a; b, bc/a,bd/a,be/a,bf/a,bg/a,bh/a;q).
$$
Here $b$ is a distinguished parameter which can be interchanged with
$ c,d,e,f,g $  or $h$ to give different $\Phi$'s. We use the notation
$$ \Phi \ = \ \Phi(a;b,c,d,e,f,g,h;q) $$
to denote any one of these seven possible complementry pairs.

\beginsection 2. Contiguous relation.

It has been show that, irrespective of which is the distinguished
paramster, $\Phi$ satisfies the three term contiguous relation [5],[6]
$$
\leqalignno{
& { g(1-h) (1-{a\over h}) 
(1- {aq\over h})(1-{aq\over gb}) 
(1- {aq\over gc})
(1- {aq\over gd})
(1- {aq\over ge})
(1- {aq\over gf})\over
(1- {hq\over g}) } \cr
&\hskip  2in \times \big[\Phi (g -, h+) -\Phi \big] \cr
- & { h(1-g) (1-{a\over g}) 
(1- {aq\over g}) 
(1- {aq\over hb})
(1- {aq\over hc})
(1- {aq\over hd})
(1- {aq\over he})
(1- {aq\over hf})\over
(1- {gq\over h}) } \cr
&\hskip  2in \times \big[\Phi (h -, g+) -\Phi \big] \cr
- & {aq\over h} (1-{h\over g}) (1-{gh\over aq})
(1-b) (1-c) (1-d) (1-e) (1-f) \Phi \ = \ 0 \ .\cr}
$$
Note that we are using the notation $\Phi (g -, h +)$ to mean a
$\Phi$ with $g$ replaced by $g/q$ and $h$ replaced by
$hq$.

The proof of this contiguous relation is based on extending the
contiguous relations derived in [5] for a terminating $\phi$ to
nonterminating complementary pairs $\Phi$. This is made possible
by using the transformation in exercise 2.30 of Gasper and Rahman [2]
to replace the transformation of exercise 2.19 used in [5].

\beginsection 3. Difference equation solutions.

In the contiguous relation of section 2 we put $ h = q^{-n}$ and 
$ g = sq^{n-1} $, where
$$ s \ := \ a^3q^3/bcdef $$
takes care of the balance condition. 

\noindent We now renormalize and reexpress the contiguous relation as
$$ X_{n+1} -a_n X_n +b_n X_{n-1} \ = \ 0 ,$$
$$ a_n \ = \ A_n+B_n+{sq^{2n}\over aq}{(1-{s\over aq^2})(1-b)(1-c)
(1-d)(1-e)(1-f)\over (1-aq^{n+1})(1-{s\over a}q^{n-2})},
$$
$ b_n \ = \ A_{n-1}B_n $   \quad with
$$ A_n \ = \ { (1-sq^{n-1})(1-{s\over aq}q^{n})(1-{a\over b}q^{n+1})
(1-{a\over c}q^{n+1}) \cdots (1-{a\over f}q^{n+1})\over
(1-sq^{2n})(1-s^{2n-1})(1-aq^{n+1})}
$$
$$
B_n \ = \ {q(1-q^n)(1-aq^n)(1-{bs\over a}q^{n-2})(1-{cs\over a}q^{n-2})
\cdots (1-{fs\over a}q^{n-2})\over
(1-sq^{2n-1})(1-sq^{2n-2})(1-{s\over a}q^{n-2})}
$$
and having the solution
$$
X_n^{(1)} \ = \  {(sq^{2n-1}, aq^{n+1})_\infty \over
(sq^{n-1},{s\over a}q^{n-1},{a\over b}q^{n+1},{a\over c}q^{n+1},\cdots,
{a\over f}q^{n+1})_\infty} $$
$$ \times \phi(a; b,c,d,e,f, sq^{n-1}, q^{-n} ; q) . $$

There is a remarkable hidden symmetry that yields additional solutions [5].
If we make the replacements
$$
(a,b,c,d,e,f,sq^{n-1},q^{-n}) \to (q/a,q/b,q/c,q/d,q/e,q/f,q^{2-n}/s,
q^{n+1}) $$
then
$$ a_n \to a_n \quad, \quad b_n \to b_{n+1} .$$
Applying this to the solution $X_n^{(1)}$ yields solutions 
$$ 
X_n^{(p)} \ = \ {(sq^{2n-1},sq^n/a)_\infty \over
(q^{n+1},aq^n,{bs\over a}q^{n-1},{cs\over a}q^{n-1}, \cdots,{fs\over a}q^{n-1})_\infty}
$$
$$ \times \Phi^{(p)}\big({q\over a};{q\over b},{q\over c},{q\over d},
{q\over e},{q\over f},{q^{2-n}\over s},q^{n+1};q) , $$
$ \quad \quad p \ = \ q/b,q/c,q/d,q/e,q/f,q^{2-n}/s $ or $ q^{n+1}$.

\noindent Because of its simpler large n asymptotics we choose
$$ X_n^{(2)} \quad = \quad X_n^{(q^{n+1})} $$
as a second solution.

\beginsection 4. Asymptotics and Pincherle's theorem. 

We now examine the large $n$ asymptotics of the difference equation.
Since we are assuming $ |q| \ < \ 1 $, it is easy to see that
$$ \lim_{n\to \infty} a_n \ = \ 1 \quad,\lim_{n\to \infty} b_n \ = \ q.
$$
It follows that 
$$ X_{n+1} -a_n X_n + b_n X_{n-1} \ = \ 0 $$
has solutions with large n asymptotics
$$ X_n \ \approx \ const \quad \hbox{or} \quad X_n \ \approx \ const q^n .$$
Clearly, the latter yields minimal (subdominant) asymptotics. We are
particularly interested in such a minimal solution which we denote by 
$ X_n^{(min)} $ because of the 
\proclaim Theorem (Pincherle, 1894, see [3],[8]).
$$ {1\over\displaystyle a_0 -{b_1\over \displaystyle a_1 -{b_2\over\displaystyle \ddots}}} ={X_0^{(min)}\over
a_0X_0^{(min)}-X_1^{(min)}}.$$
\par

\beginsection 5. The minimal solution and continued fraction.

We now construct the minimal solution by examining the large $n$ 
asymptotics of the solutions $X_n^{(1)}$ and $X_n^{(2)}$.
 
For $X_n^{(1)}$ it is easy to see that
$$ \lim_{n\to \infty} X_n^{(1)}=W_1 $$
where
$$ W_1 \ =\ {}_8\phi_7\left( \matrix {a,&q\sqrt{a},&-q\sqrt{a},&b,&c,
&d,&e,&f \cr{}&\sqrt{a},&-\sqrt{a},&{aq\over b},&{aq\over c},&{aq\over d},&{aq\over e},&{aq\over f},\cr}
; { a^2q^2\over bcdef}\right).$$
For $X_n^{(2)}$ it is not immediate but a calculation yields 
$$ \lim_{n\to \infty} X_n^{(2)}=W_2 $$
where
$$ W_2 \ =\ {}_8\phi_7\left( \matrix {{q\over a},&q\sqrt{{q\over a}},
&-q\sqrt{{q\over a}},&{q\over b},&{q\over c},&{q\over d},&{q\over e},
&{q\over f} \cr {}&\sqrt{{q\over a}},&-\sqrt{{q\over a}},&{bq\over a},
&{cq\over a},&{dq\over a},&{eq\over a},&{fq\over a} \cr}; {bcdef\over 
a^2q}\right).$$

 A minimal solution is therefore given by
 $$ X_n^{(min)} \ = \ W_2X_n^{(1)}-W_1X_n^{(2)} $$
 and we may now apply Pincherle's theorem. The result is
\proclaim Theorem 1.
$$ 
\leqalignno{
&{1\over\displaystyle a_0 -{b_1\over \displaystyle a_1 -
{b_2\over\displaystyle \ddots}}}=
{(1-s/q)(1-a/q)\over q(1-{s\over aq})(1-{a\over b})
(1-{a\over c})\cdots (1-{a\over f})}\cr
&\times\Bigl\lbrack \phi({q\over a};q,{q^2\over s},{q\over b},{q\over c},
\cdots,
{q\over f};q)+\Pi_1 \phi(aq;q,{aq^2\over s},{aq\over b},
\cdots,{aq\over f};q)\cr
&\quad \quad-{(q,a,aq,{bs\over aq},\cdots,{fs\over aq})_\infty 
\over({s\over q},{s\over a}, {s\over aq},{aq\over b},\cdots,{aq\over f})
_\infty}{W_1\over W_2}\Bigr\rbrack/(1+\Pi_0),\cr
&\Pi_n:={(q^2/a,q/b,\cdots,q/ f,q^{3-n}/ s,
aq^{n-1},sq^{2n-2},bq^n,\cdots,fq^n)_\infty
\over(aq^{2n},q^{1-n} /a,bq /a,\cdots,fq/ a,aq^2/ s,sq^{n-1}/ a,aq^n/ b
,\cdots,aq^n/ f)_\infty}.\cr}$$
\par
However, a remarkable simplification occurs if the fraction terminates.
\proclaim Corollary 2.
If one of ${aq\over b},{aq\over c},\cdots,{aq\over f}=q^{-N},N=0,1,
\cdots$, then
$$ 
\leqalignno{
&{1\over\displaystyle a_0 -{b_1\over \displaystyle a_1 -
{b_2\over\displaystyle \ddots}}}={aq(1-aq)\over s(1-b)(1-c)(1-d)(1-e)(1-f)}\cr
&\times \phi(aq;q,aq^2/s,aq/b,aq/c,aq/d,aq/e,aq/f;q).\cr}
$$
\par

\beginsection 6. Two limit cases.

 In corollary 2 we can take the limit $N\to \infty$ to obtain two
 limit cases at the ${}_8\phi_7$ level.
 With $f=aq^{N+1}$ , $s=a^2q^{2-N}/bcde$ and $N\to \infty$ in corollary 2 we have
 \proclaim Corollary 3.
$$ 
\leqalignno{
&{1\over\displaystyle c_0 -{d_1\over \displaystyle c_1 -
{d_2\over\displaystyle \ddots}}}={bcde\over a^2q^2}{(1-aq)\over 
(1-b)(1-c)(1-d)(1-e)}\cr
&\times{}_8\phi_7 \left(\matrix{aq,&q\sqrt{aq},&-q\sqrt{aq},&q,&{aq\over b},\cdots,&{aq\over e}\cr
&\sqrt{aq},&-\sqrt{aq},&aq,&bq,\cdots,&eq\cr};{bcde\over a^2q}\right)\cr}
$$
where
$$c_n=-{(1-{aq^{n+1}\over b})(1-{aq^{n+1}\over c})(1-{aq^{n+1}\over d})
(1-{aq^{n+1}\over e})\over (1-aq^{n+1})}$$
$$
-q(1-q^n)(1-aq^n)(1-{a^2q^{n+1}\over bcde})
+{a^2q^{2n+2}\over bcde}{(1-b)(1-c)(1-d)(1-e)\over (1-aq^{n+1})},
$$
$$
d_n=q(1-q^n)(1-{aq^n\over b})(1-{aq^n\over c})(1-{aq^n\over d})
(1-{aq^n\over e})(1-{a^2q^{n+1}\over bcde}).
$$
\par
If $f=aq^{N+1}$ , $e=a^2q^{2-N}/bcds$ (with $s$ now an independent
parameter) and $N\to \infty$ in corollary 2, then we recover the
continued fraction result associated with Askey-Wilson polynomials [4],[7].
\beginsection 7. Remarks.

\indent 1. Theorem 1 generalizes a result in [5] which considered only $s=q^n,
\quad n=1,2,\cdots$.

2. The case $s=q^2$ in corollary 2 is the q-analogue of Ramanujan's
famous Entry 40 in Chapter 12 of his Notebook 2 (see [5]).

3. A previous evaluation of the continued fraction in corollary 3 was given
in [5]. It was was much more complicated looking.

4. The continued fractions in corollaries 2 and 3 are not associated
with orthogonal polynomials but they do yield results for biorthogonal
rational functions. This will be reported on elsewhere.

5. Theorem 1 can be generalized slightly. In section 3 we put $h=q^{-n}$.
One could consider, instead, the associated case with $h=\epsilon q^{-n}$. 
However, the formulas become considerably more complicated. 
\bigskip 

\centerline {\bf References }
\baselineskip=12pt
\frenchspacing
\vglue .25in

\item{1. }
R. Askey and J. Wilson,
Some basic hypergeometric orthogonal polynomials that generalize
Jacobi polynomials,
{\it Memoirs Amer. Math. Soc.} {\bf 319} (1985) 1--55.
\medskip

\item{2. } G. Gasper and M. Rahman,
{\it Basic Hypergeometric Series,}
Cambridge Univ. Press, Cambridge, 1990.
\medskip

\item{3. }
W. Gautschi,
Computational aspects of three-term recurrence relations,
{\it SIAM Rev.} {\bf 9} (1967), 24--82.
\medskip

\item{4. } D.P. Gupta and D.R. Masson,
Exceptional  $ q $-Askey-Wilson polynomials and continued fractions,
{\it Proc. A.M.S.} {\bf 112} (1991), 717--727.
\medskip

\item{5. } 
D.P. Gupta and D.R. Masson,
Watson's basic analogue of Ramanujan's Entry 40 and its
generalization,
{\it SIAM J. Math. Anal.}, to appear.
\medskip

\item{6. }
D.P. Gupta and D. R. Masson, 
Contiguous relations and continued fractions,
{} manuscript in preperation.
\medskip

\item{7. } M.E.H. Ismail and M. Rahman,
Associated Askey-Wilson polynomials,
{\it Trans. Amer. Math. Soc.} {\bf 328} (1991), 201--239.
\medskip

\item{8. }
W.B. Jones and W.J. Thron, {\it Continued Fractions: Analytic
Theory and Applications,} 
Addison-Wesley, Reading, Mass., 1980.
\medskip

\end
{(1-c)  (1-{a\over c})
(1-{aq\over c})
(1-{aq\over bd})
(1-{aq\over be})
(1-{aq\over bf})\over 
(1- {cq\over b}) } \ 
\big[ \ephis (b-, c+) - \ephis \big]  \cr-
& {bcdef\over a^2 q} 
{(1-b)  (1-{a\over b})
(1-{aq\over b})
(1-{aq\over cd})
(1-{aq\over ce})
(1-{aq\over cf})\over 
(1- {bq\over c}) } \ 
\big[ \ephis (c-, b+) - \ephis \big]  \cr
& \hskip .5truein -
 {aq\over c} (1-{c\over b}) (1-{bc\over aq})
(1-d) (1-e) (1-f) \ephis \ = \ 0 \ ,\cr}
$$
where
$$
\leqalignno{
\ephis \ = \ & \ephis 
\left (  
 { {a,q\sqrt a, -q\sqrt a, b, c,d,e,f} 
   \atop 
 {\displaystyle \sqrt a, -\sqrt a, {aq\over b}, {aq\over c},
 {aq\over d}, {aq\over e},
 {aq\over f}}} \ ; \  
  {a^2 q^2\over bcdef} \right)  \cr
  = \ & W(a; \ b, c, d, e, f ) \quad \hbox{say} \ . \cr}
$$
For non-terminating  $ \ephis $,  the convergence condition is
$ | a^2 q^2/ bcdef | < 1 $.
In (2.3) we now put
$$
\leqalignno{
a \ = \ &  {\alpha \beta \gamma u\over q}, \quad
b \ = \ q^{-n}/ \epsilon\ , \quad
c \ = \ \epsilon s q^{n-1}  \ , \quad
d \ = \ \alpha u \ , \cr
e \ = \ & \beta u \ , \quad
f \ = \ \gamma u ,
s \ = \  \alpha\beta\gamma \delta \ .\cr}
$$
Renormalizing and simplifying, we obtain the associated Askey-Wilson
recurrence
$$
X_{n+1}^{(1)}  - (z - a'_n ) X_n^{(1)}  + b_n^{\prime 2}  X_{n-1}^{(1)}  \ = \ 0
\leqno(2.5)
$$
with
$$
\eqalign{
z \ = \ & {u+ u^{-1}\over 2 } \ , \cr
a'_n \ = \ & -A'_n - B'_n + {\alpha\over 2} + {1\over 2\alpha} 
 \cr
b^{\prime 2}_n \ = \ & A'_{n-1}  B'_n  
 \cr
A'_n \ = \ &
{ (1-sq^{n-1} \epsilon ) 
  (1-\alpha\beta \epsilon q^n ) 
  (1-\alpha\gamma \epsilon q^n ) 
  (1-\alpha\delta \epsilon q^n ) \over
  2\alpha (1-sq^{2n-1} \epsilon^2) (1-sq^{2n}\epsilon^2) } \ , \cr
B'_n \ = \ &
{ \alpha (1-\epsilon q^n ) 
  (1-\beta \gamma \epsilon q^{n-1} ) 
  (1-\beta\delta \epsilon q^{n-1} ) 
  (1-\gamma\delta \epsilon q^{n-1} )  \over
  2 (1-s\epsilon^2q^{2n-2} ) (1-s\epsilon^2q^{2n-1}) }\  \cr}
$$
and a solution of (2.5) is therefore
$$
\leqalignno{
X_n^{(1)} (u) \ = \ &
 \Big( {u\over 2} \Big)^n\ 
{ \big(su\epsilon q^n/\delta\big)_\infty
 \big(s\epsilon^2q^{2n-1}\big)_\infty
\over
 \big(s\epsilon q^{n-1}\big)_\infty
 \big(\delta\epsilon q^n/u\big)_\infty
\big(\alpha \beta\epsilon q^n , 
\alpha \gamma \epsilon q^n , 
\beta \gamma \epsilon q^n \big)_\infty }\cr
& \ \times W \Big( {\alpha\beta \gamma u\over q}\ ; \ 
{q^{-n}\over \epsilon} , \epsilon sq^{n-1} , \alpha u, \beta u, 
\gamma u \Big) \ , \cr}
$$
with convergence condition   
$ | q/\delta u| < 1 $.

The solution  $ X_n^{(1)} $  is a generalization of the 
Askey-Wilson solution (1.1) to the  $\epsilon = 1 $  case of
(2.5).  This can be verified through an application of Watson's
formula ([5], III.18, p.242) connecting a terminating $\fphit $  with
an $ \ephis $.

A second linearly independent solution of (2.5) may be obtained with 
the help of the reflection symmetry transformation 
$ v \to -v -1 $,  
$ (\alpha, \beta, \gamma, \delta) \to 
( q /\alpha,
q /\beta,
q /\gamma,
q /\delta)$.
Using (see also [10])  
$$
\eqalign{
b_{-v-1}^2 
( q /\alpha,
q /\beta,
q /\gamma,
q /\delta)\ = \ & 
b_{ v+1}^2 (\alpha, \beta, \gamma, \delta) \cr
a_{-v-1}
( q /\alpha,
q /\beta,
q /\gamma,
q /\delta)\ = \ & 
a_v (\alpha, \beta, \gamma, \delta) \cr}
$$
and renormalizing, we arrive at the solution
$$
\leqalignno{
X_n^{(2)} (u) \ = \ &
\Big({ u\over 2}\Big)^n\ 
{ \big(s\epsilon^2 q^{2n-1}\big)_\infty
 \big(\epsilon \delta u q^{n+1}\big)_\infty
\over
 \big(\epsilon q^{n+1}\big)_\infty
 \big(\beta\delta\epsilon q^n ,
\gamma \delta\epsilon q^n , 
\alpha \delta\epsilon q^n \big)_\infty 
\big( {s\epsilon\over \delta u} q^{n-1} \big)_\infty } \cr
& \ \times W \Big( {q^2u \over \alpha \beta \gamma }\ ; \ 
 \epsilon q^{n+1} ,{q^{-n+2}\over \epsilon s } , {qu \over\alpha } ,
 {qu \over\beta} , 
{qu \over\gamma}  \Big) \ , \cr}
$$
with convergence condition  $ |\delta /u|< 1 $.

In order to obtain a third solution of (2.5) we make the
following alterations in (2.2):
\noindent
Interchange  $ c $  and  $ h $,  replace  $ a,b,c,\ldots , g $
by  $ h^2/ a $,  $ hb /a $,
$ hc/ a , \ldots , hg /a $  respectively, put
$ h = q^{-n} $,  $ n $  being a positive integer, and reverse
the  $ \tphin $  series.  This gives
$$
\leqalignno{
& { {b\over a} (1-h)(1-{h\over a}) 
(1-{hq\over a})
(1-{aq\over bd})
(1-{aq\over be})
(1-{aq\over bf}) 
(1-{aq\over bg}) 
(1-{aq\over bc})\over
(1- {aq\over b}) } \ 
  \cr
&\hskip  1in \times \big[c_2 \phi_+ (b -) - c_1 \phi\big]\cr
- & {( 1- {bh\over a} )
(1-{h\over b}) 
(1-{hq\over b})
(1-{q\over d})
(1-{q\over e})
(1-{q\over f}) 
(1-{q\over g}) 
(1-{q\over c})\over
(1- {bq\over a}) } \cr 
&\hskip  1in \times \big[c_3 \phi_- (b +) - c_1 \phi\big]\cr
- &{q\over a}
( 1- {a\over b} )
(1-{b\over q}) 
(1-{hd\over a})
(1-{he\over a})
(1-{hf\over a})
(1-{hg\over a}) 
(1-{hc\over a}) 
c_1 \phi \ = \ 0 \ , \cr}
$$
where
$$
\leqalignno{
& c_1 \ = \ 
{ (\sqrt a )_n (-\sqrt a)_n 
\big( {aq\over b}\big)_n
\big( {aq\over c}\big)_n \ldots
\big( {aq\over g}\big)_n
\big( aq^{n+1} )_n\over
(a)_n (q\sqrt a)_n (-q\sqrt a)_n (b)_n (c)_n \ldots (g)_n } \ 
(-1)^n q^{n(n-1)/2} \cr
& c_2 \ = \ (-1)^{n-1}q^{(n-1)(n-2)/2} \cr
&\times { (q\sqrt a )_{n-1} (-q\sqrt a)_{n-1} 
\big( {aq^3\over b}\big)_{n-1}
\big( {aq^2\over c}\big)_{n-1} \ldots
\big( {aq^2\over g}\big)_{n-1}
\big( aq^{n+2}\big )_{n-1}\over
(aq^2)_{n-1} (q^2\sqrt a)_{n-1} (-q^2 \sqrt a)_{n-1} (b)_{n-1} 
(cq)_{n-1} \ldots (gq)_{n-1} } \cr 
 & c_3 \ = \ (-1)^{n+1}q^{n(n+1)/2}\cr
& \times{ \big({\sqrt a \over q}\big)_{n+1} \big(-{\sqrt a\over q}\big)_{n+1} 
\big( {a\over bq}\big)_{n+1}
\big( {a\over c}\big)_{n+1} \ldots
\big( {a\over g}\big)_{n+1}
\big( aq^n \big)_{n+1}\over
\big({a\over q^2}\big)_{n+1} (\sqrt a)_{n+1} (-\sqrt a)_{n+1}
(b)_{n+1} 
\big({c\over q}\big )_{n+1} \ldots \big({g\over q}\big)_{n+1} }  
 . \cr}
$$
 Note that $\phi_+(b-) $ and $\phi_-(b+) $ refer to $ \phi $ with $ (a, c, d, 
e, f, g, h) $ replaced by $(aq^2, cq, dq, eq, fq, gq, hq) $ and $(aq^
{-2}, c/q, d/q, e/q, f/q, g/q, h/q) $ respectively.

Taking the limit  $ h = q^{-n}\to \infty $,
$ g = a^3 q^{2+n} /bcdef \to 0 $,  we obtain
$$
\leqalignno{
& 
{b(1-{aq\over bd})
(1-{aq\over be})
(1-{aq\over bf})
(1-{aq\over bc})\over
a(1-{aq\over b})} \cr
\times & \left[
{a^2 q^2\over bcdef}\ 
{(1-aq) (1-aq^2) (1-c) (1-d) (1-e) (1-f) \over
(1-{aq\over b})
(1-{aq^2\over b})
(1-{aq\over c})
(1-{aq\over d})
(1-{aq\over e})
(1-{aq\over f})} \ 
W_+ (b-) - W \right ]\cr
& \hskip .5truein -{
(1-{q\over d})
(1-{q\over e})
(1-{q\over f}) 
(1-{q\over c})\over
(1- {bq\over a})} \cr 
\times &  \left[ {bcdef\over a^2 q^2} 
{ (1-{a\over bq})
(1-{a\over b})
(1-{a\over c})
(1-{a\over d})
(1-{a\over e})
(1-{a\over f})\over 
(1- {a\over q})  (1-a) 
(1-{c\over q})
(1-{d\over q})
(1-{e\over q})
(1-{f\over q})} \ 
W_- (b+) - W \right ]\cr
& \hskip .5truein -{q\over a} 
(1-{a\over b})
(1-{b\over q})
\big(1-{a^2 q^2\over bcdef}\big) W \ = \ 0 . \cr}
$$

In (2.9), writing
$$
\leqalignno{
a \ = \ & {\beta \gamma \delta \epsilon^2\over u} q^{2n} , \quad
b \ = \  {q\over \alpha u} , \quad
c \ = \  q^{n+1}\epsilon  , \quad
d \ = \  \beta \delta \epsilon q^n   \cr
e \ = \ &  \gamma \delta \epsilon  q^n, \quad
f \ = \  \beta \gamma \epsilon q^n  , \quad
s \ = \  \alpha \beta \gamma \delta \cr}
$$
and renormalizing, we obtain a third solution to (2.5) given by
$$
\leqalignno{
& \quad X_n^{(3)} (u)\ = \ 
{1\over (2u)^n}{\big(s\epsilon^2 q^{2n}\big)_\infty\big(s\epsilon^2 
q^{2n-1}\big)_\infty \over\big(\epsilon q^{n+1}\big)_\infty
\big(\epsilon sq^{n-1}\big)_\infty}
\cr 
& \qquad \qquad\times {
\big( {\beta\gamma \delta\epsilon\over u} q^n \big)_\infty 
\big( {\epsilon \beta q^{n+1}\over u},
 {\epsilon \gamma q^{n+1}\over u},
 {\epsilon \delta q^{n+1}\over u}\big)_\infty
 \over
\big( \alpha\beta \epsilon q^n , 
\alpha \gamma\epsilon q^n,
\alpha \delta\epsilon q^n,
\beta \gamma\epsilon q^n,
\beta \delta\epsilon q^n,
\gamma \delta\epsilon q^n \big)_\infty 
\big( {\beta \gamma \delta \epsilon^2\over u} q^{2n+1} \big)_\infty }
 \cr
&\qquad \qquad \qquad  \times W \Big( {\beta \gamma \delta \epsilon^2\over u}
q^{2n} \ ; \ 
{q\over \alpha u} , 
\epsilon q^{n+1} , 
\beta \delta\epsilon q^n , 
\gamma \delta \epsilon q^n ,
\beta \gamma \epsilon q^n  \Big) \ , \cr}
$$
with convergence condition  $ | \alpha/ u | < 1 $.
A transformation ([5], III.23, p.243) may be applied to  $ X_n^{(3)}$
to obtain this solution in a form which is symmetric in the
parameters  $ \alpha, \beta, \gamma, \delta $.  We get the solution
(omitting constant factors)
$$
\leqalignno{
& X_n^{(4)} (u)
= \ 
{1\over (2u)^n}\cr 
 \times & {\big(\epsilon ^2 sq^{2n-1}\big)_\infty
\big( {\alpha\epsilon\over u} q^{n+1} ,
 {\beta \epsilon \over u}  q^{n+1},
 {\gamma \epsilon \over u}  q^{n+1},
 {\delta \epsilon \over u}  q^{n+1}\big)_\infty
 \over
 \big(\epsilon q^{n+1}\big)_\infty
 \big({q^{n+2}\epsilon\over u^2}\big)_\infty
\big( \alpha\beta \epsilon q^n , 
\alpha \delta\epsilon q^n,
\alpha \gamma\epsilon q^n,
\beta \gamma\epsilon q^n,
\beta \delta\epsilon q^n,
\gamma \delta\epsilon q^n \big)_\infty }
 \cr
&\ \times W \Big( 
{q^{n+1}\epsilon \over u^2} ; \ 
\epsilon q^{n+1},
{q\over \alpha u} , 
{q\over \beta u},
{q\over \gamma u},
{q\over \delta u} \Big)  \cr}
$$
where the convergence condition is  $ |sq^{n-1}\epsilon | < 1 $.

The large  $ n $ behaviour of  $ X_n^{(4)} (u) $  is easily seen
to be
$$
\leqalignno{
& X_n^{(4)} (u)\ 
\ttopsim{n\to\infty}
{1\over (2u)^n}\cr
\noalign{\hbox{while the corresponding solution}}
& X_n^{(4)} (1/u)
\ \sim \ \big({u\over 2}\big)^n . \cr}
$$
It is evident that  $ X_n^{(4)} (u) $  is a sub-dominant or
minimal solution
for  $ |u| > 1 $,  $ |sq^{n-1} \epsilon | < 1 $   while 
$ X_n^{(4)} (1 /u) $ 
is sub-dominant (minimal) for $ | u | < 1 $,
$ |sq^{n-1} \epsilon | < 1 $.

We proceed to get another solution by making the following parameter
replacements in (2.9):
$$
\eqalign{
a \ = \ & {q^{-2n+1}\over \beta\gamma \delta u\epsilon^2}, \quad
b \ = \  {\alpha \over u},\quad
c \ = \  {q^{-n}\over \epsilon } ,\quad
d \ = \  {1\over \beta\delta \epsilon } \ q^{-{n+1}},\cr
e \ = \  & {1\over \gamma\delta \epsilon } \ q^{-{n+1}},\quad
f \ = \   {1\over \beta \gamma \epsilon } \ q^{-{n+1}}\ .\cr}
\leqno(2.14)
$$
Note that this is just the reflection transformation that was 
used to obtain the solution $ X_n^{(2)} $ from $ X_n^{(1)} $ but
now it is being applied to $ X_n^{(3)} $.

After renormalization, omitting constant factors
and using the transformation [5, III.23, p.243] we arrive at the
following solution 
$$
\leqalignno{
 X_n^{(5)} (u) 
\ = \ & {1\over (2u)^n}\ 
{\big(sq^{2n-1}\epsilon^2\big)_\infty
\big( u^2 \epsilon q^n\big)_\infty
\over 
\big(sq^{n-1}\epsilon\big)_\infty
\big(\alpha u \epsilon q^n , 
\beta u\epsilon q^n,
\gamma u \epsilon q^n,
\delta u \epsilon q^n\big)_\infty } \cr
&\quad \times W \Big( 
{q^{-n}\over \epsilon  u^2} ; \ 
{q^{-n}\over \epsilon },
{\alpha \over u} , 
{\beta\over u},
{\gamma\over u},
{\delta \over u} \Big)  \cr}
$$
where the convergence condition is
$ \big| q^{-n+2} /s\epsilon \big| < 1 $.

Note that, when $ \epsilon  =1 $, the solution (2.15) is proportional
to the Askey-Wilson polynomial $ P_n(z; \alpha ,\beta , \gamma , \delta
) $ in (1.1) and can be used to give $ P_n(z) $ in a form which is
explicitly symmetric in its parameters. That is
$$
P_n(z) \ = \ (2u)^{-n}{(\alpha u, \beta u, \gamma u, \delta u, s/q)_n
\over (s/q)_{2n}(u^2)_n}W(q^{-n}/u^2;q^{-n},{\alpha \over u},
{\beta \over u}, {\gamma \over u}, {\delta \over u}) .
$$

A solution in which the convergence condition is independent
of  $ z $  and also independent of the parameters
$ \alpha, \beta, \gamma, \delta $  can be obtained from (2.12)
by applying the transformation
([5], III.24, p.243).  This gives
$$
\leqalignno{
  & X_n^{(6)} (u) 
\ = \  {1\over (2u)^n}\cr 
\times & {\big(s\epsilon^2 q^{2n-1}\big)_\infty \big(
s\epsilon q^n/ \alpha u,
s\epsilon q^n/ \beta  u,
s\epsilon q^n/ \gamma  u,
s\epsilon q^n/ \delta  u\big)_\infty\over
\big(\epsilon sq^{n-1}\big)_\infty
\big(s\epsilon  q^n/u^2 \big)_\infty
\big( \alpha \beta\epsilon q^n , 
\alpha\gamma \epsilon q^n,
\alpha\delta \epsilon q^n,
\beta\gamma \epsilon q^n,
\beta\delta \epsilon q^n,
\gamma\delta \epsilon q^n\big)_\infty} \cr
&\quad \times W \Big( 
{s\epsilon q^{n-1}\over u^2} ; \ 
{\epsilon s q^{n-1}},
{\alpha \over u} , 
{\beta\over  u},
{\gamma\over  u},
{\delta \over u} \Big)  \cr}
$$
with convergence condition  $ \big| q^{n+1}\epsilon \big| < 1 $.

Since  $ X_n^{(6)} (u) \topsim{n\to\infty} 
{1\over (2u)^n} $,
 we have that $  X_n^{(6)} (u) $  is a sub-dominant(minimal) solution for
$ |u | > 1 $  while
  $  X_n^{(6)} (1/u) $  is a sub-dominant(minimal) solution for 
$ |u | < 1 $.
 Note that  $ X_n^{(6)} $  is really the same solution as
$ X_n^{(4)} $,  but analytically continued.

Summarizing we have

\proclaim Theorem 1.  The associated Askey-Wilson equation (2.5)
has solutions  $ X_n^{(i)} (u) $ and $ X_n^{(i)}(1/u) $ ,
$ i = 1,2,\ldots , 6 $  given by
(2.6), (2.7), (2.11), (2.12), (2.15) and (2.16).  If  $ |u| > 1 $,  
then  
$ X_n^{(3)} (u) $,   
$ X_n^{(4)} (u) $  and  
$ X_n^{(6)} (u) $ 
(which are connected by  general  $ \ephis $  transformations) each
represent a minimal solution
$ X_n^{(s)} (z) $ 
for the parameter values  
$ \big|{\alpha\over u}\big| < 1 $,
$ |\epsilon sq^{n-1}| < 1 $,
$ |\epsilon q^{n+1}| < 1 $
respectively.\par
\vglue .5truein

\beginsection 3. Continued fraction representations.

The continued fraction associated with the difference equation
(2.5) is
$$
CF (z) \ :=\ z - a'_0 + \bigk\limits_{n=1}^\infty
\left( {-b_n^{\prime 2} \over z -a'_n } \right)  \ .
\leqno(3.1)
$$
If  $ b_n^{\prime 2} \not= 0 $,  $ n \ge 1 $,  by Pincherle's
theorem [6], [17] and Theorem 1 of Section 2 we have
$$
{1\over CF (z)} \ =\ 
{ X_0^{(s)} (z)\over
b_0^{\prime 2} X_{-1}^{(s)} (z)} \ .
\leqno(3.2)
$$
Using  
$ X_0^{(s)} (z) = X_0^{(4)} (u)$,
$ X_{-1}^{(s)} (z) = X_{-1}^{(4)} (u)$,
$ |u| > 1, z={ u + u^{-1} \over 2} $,
we have from (2.12),
$$
\leqalignno{ 
{1\over CF (z)} \ =\ 
& {2\over u} 
{ (1-\epsilon q/ u^2) 
(1- s\epsilon^2/ q^2)
(1- s\epsilon^2/ q)
\over
( 1-\alpha\epsilon / u)
( 1-\beta\epsilon/ u)
( 1-\gamma\epsilon/ u)
( 1-\delta\epsilon/ u)
( 1-s\epsilon/ q^2) }\cr
& \times 
{W \Big( 
q\epsilon/ u^2 ; \ q\epsilon,
q/ \alpha u,
q/ \beta u,
q/ \gamma  u,
q/ \delta  u\Big)
\over
W \Big( 
\epsilon/ u^2 ; \ \epsilon,
q/ \alpha u,
q/ \beta u,
q/ \gamma  u,
q/ \delta  u\Big)} & (3.3)\cr}
$$
for  $ |u| > 1 $  and  $ |s\epsilon/ q^2 | < 1 $.
For $|u| < 1, \,|s\epsilon/q^2| < 1 $ we replace $u$ in
(3.3) by $ 1/u $.

 Other representations for the continued fraction for 
different parameter ranges are
obtained by taking a different representation for the minimal
solution.
 For example, by taking  
$ X_0^{(s)} (z) = X_0^{(6)} (u)$ and
$ X_{-1}^{(s)} (z) = X_{-1}^{(6)} (u)$, 
$ | u | > 1 $,  we obtain the continued fraction representation
$$
\leqalignno{ 
{1\over CF (z)} \ =\ 
& {2\over u}{ 
(1- s\epsilon^2 /q^2)
(1- s\epsilon^2 /q)(1-s\epsilon/u^2q) \over
(1-\epsilon)(1-s\epsilon/\alpha uq)(1-s\epsilon/\beta uq)(1-s\epsilon
/\gamma uq)(1-s\epsilon/\delta uq)}\cr
& \times
{W \Big( 
\epsilon s/u^2q  ; \ \epsilon s/q,
 \alpha /u,
 \beta/ u,
 \gamma / u,
 \delta / u\Big)
\over
W \Big( 
s\epsilon /u^2q^2 ; \ \epsilon s /q^2 ,
\alpha /u,
\beta /u,
\gamma / u,
\delta / u\Big)} , \cr}
$$
for  $ |u| > 1 $,  $ | \epsilon |  < 1 $.
\vglue .5truein

\beginsection 4. Weight function.

We have seen above that a minimal solution exists for
$ | u | > 1 $  and  $ | u | < 1 $.  There is no minimal solution
for  $ z \in [-1,1] $,  $(|u|=1)$,  where we have a continuous
spectrum for the associated tridiagonal Jacobi matrix  $ J $
with diagonal 
$ (a'_0, a'_1, \ldots) $  and  
$ (b'_1, b'_2, \ldots) $  above and below the diagonal.
For a probability measure  $ d\omega (x;\epsilon) $
we have, in the case of real orthogonality $ ( b_{n+1}^{\prime 2} 
> 0 , a'_n\quad\hbox{real},\quad n\ge 0 )$ for the associated monic
$ q $-Askey-Wilson polynomials
$ P_n (x; \epsilon)$,
$$
\int_\R \ P_n (x;\epsilon) P_m (x;\epsilon)
d\omega (x;\epsilon) \ = \ 
\delta_{nm}  \prod\limits_{k=1}^n \ b_k^{\prime 2} \ .
\leqno(4.1)
$$
From [21], we have the representation
$$
{1\over CF(z)} \ = \ 
\int_\R \ 
{d\omega (x;\epsilon)\over z-x} \ = \ 
{X_0^{(s)} (z)\over b_0^{\prime 2} X_{-1}^{(s)} (z)}
\leqno(4.2)
$$
and for the absolutely continuous part for  $ x\in [-1,1] $,
$$
\leqalignno{
{d\omega (x;\epsilon)\over dx} \ = \ 
& {1\over 2\pi i b_0^{\prime 2} } \ 
\left( 
{X_0^{(s)} (x-i 0)\over 
X_{-1}^{(s)} (x-i 0)} -
{X_0^{(s)} (x+i 0)\over 
X_{-1}^{(s)} (x+i 0)} \right) & (4.3) \cr
= \ & {1\over 2\pi i  b_0^{\prime 2} } \ 
{\cW \big(X_{-1}^{(s)} (x+i 0),
X_{-1}^{(s)} (x-i 0) \big)\over
\big| X_{-1}^{(s)} (x-i 0)\big|^2 } , \cr}
$$
where  
$$
W (X_n, Y_n ) := X_n Y_{n+1} - X_{n+1} Y_n \ .
$$

From (2.5) we have
$$
\cW \Big(  X_{-1}^{(4)} (u), X_{-1}^{(4)} 
(1/ u)\Big) \ = \ 
\lim\limits_{n\to \infty} 
{\cW \Big( X_n^{(4)} (u), X_n^{(4)} 
(1 /u)\Big) \over
\prod\limits_{k=0}^n b_k^{\prime 2}}.
\leqno(4.4)
$$
Using (2.5) and (2.12) we then obtain
$$
\leqalignno{
\cW \Big(  X_{-1}^{(4)} (u), & X_{-1}^{(4)} 
(1 /u)\Big)  \cr
& = \  2(u-u^{-1})  
{\big({s\epsilon^2\over q^3}\big)_\infty
\big({s\epsilon^2\over q^2}\big)_\infty
\over
\Big( {\alpha \beta \epsilon \over q} ,
{ \alpha \gamma \epsilon \over q} ,
{ \alpha \delta \epsilon \over q} ,
{ \beta \gamma \epsilon \over q} ,
{ \beta \delta \epsilon \over q} ,
{ \gamma \delta\epsilon\over q} ,
{ s \epsilon\over q^2} ,
\epsilon \Big)_\infty }\ . \cr}
$$
Thus (4.3) with $$ X_n^{(s)}(x \pm i0)\ = \ X_n^{(4)}(e^{\mp i\theta}),
\quad
x\ = \ \cos\theta,\quad u \ = \ e^{i\theta} $$

 \noindent gives for  $ -1\le x \le 1 $,  
$ \big| s\epsilon /q^2 \big| < 1 $,
$$
\leqalignno{
& {d\omega (x;\epsilon )\over dx} \ = \ 
{2\sqrt {1-x^2}\over \pi }\ 
{(1 - s\epsilon^2/q)
(1 - s\epsilon^2 /q^2)^2\over
(1-s\epsilon/q^2)} & (4.6)\cr
\times \ &
{\big(
\epsilon q /u^2, \epsilon q u^2,
\alpha \beta \epsilon,
\alpha \gamma \epsilon,
\alpha \delta \epsilon,
\beta\gamma \epsilon,
\beta\delta \epsilon,
\gamma\delta \epsilon,
\epsilon q \big)_\infty\over
\big(\alpha\epsilon /u ,
\alpha\epsilon u,
\beta \epsilon /u,
\beta \epsilon u,
\gamma\epsilon /u,
\gamma\epsilon u,
\delta\epsilon /u,
\delta\epsilon u,
s\epsilon /q^2\big)_\infty}  \cr
\times \ &
{1\over \Big| W \big( 
\epsilon /u^2; \ 
q /\alpha u,
q /\beta u,
q /\gamma u,
q /\delta u,
\epsilon \big)\Big|^2} \ .\cr}
$$
This checks with the weight function obtained by Ismail and
Rahman ([16], (4.31), p.218). When  $ \epsilon = 1 $,  this
reduces to the Askey-Wilson weight function [1].
Also from the positivity of the denominators in (3.4) there
is no discrete spectrum if   $ -1 < \epsilon < 1$,
$ |\alpha |, | \beta|, |\gamma |, |\delta |< |q|^{1/2} $, [16].
Other conditions for the absence of a discrete spectrum
may be deduced from (3.3) and (3.4). For example, from (3.3),
there is no discrete spectrum if $ -|q^2/s| < \epsilon < |q^2/s| $,
$ | \alpha |, |\beta |, |\gamma |, |\delta | > |q|^{1/2} $ .

Next we use (4.3) with
$$
{d\omega (x;\epsilon )\over dx} \ = \ 
{1\over 2\pi i b_0^{\prime 2}}
\left( 
{X_0^{(6)} (u)\over 
X_{-1}^{(4)} (u)} -
{X_0^{(6)} (1/u)\over 
X_{-1}^{(4)} (1/u)} \right) \ ,\quad x \ = \ \cos\theta,\ u \ 
= \ e^{i\theta}.
$$
The right side simplifies to
$$
\leqalignno{
& {1\over 2\pi i}
{\big(1 - s\epsilon^2/ q^2\big)
\big(1 - s\epsilon^2/ q\big) ( \epsilon q)_\infty\over
\big(s\epsilon/ q^2\big)_\infty } \cr
\times \ & 
\Bigg[ 2u
{( \epsilon q u^2)_\infty
\Big( {s\epsilon\over \alpha} u ,
{s\epsilon\over \beta} u ,
{s\epsilon\over \gamma} u ,
{s\epsilon\over \delta} u \Big)_\infty
W \big({su^2 \epsilon\over q} ; \ 
{s\epsilon\over q},
\alpha u, \beta u, \gamma u, \delta u\big)\over
(su^2 \epsilon)_\infty
(\alpha u\epsilon, \beta u\epsilon, \gamma u\epsilon, 
\delta u\epsilon)_\infty
W \big( \epsilon u^2;\ \epsilon, 
{qu\over \alpha},
{qu\over \beta},
{qu\over \gamma},
{qu\over \delta}\big)} \cr
&  - {2\over u}
{\big(\epsilon q/ u^2\big)_\infty
\big( s\epsilon/ \alpha u,
s\epsilon /\beta u,
s\epsilon /\gamma u,
s\epsilon /\delta u\big)_\infty
W \big( {s\over u^2}
{\epsilon\over q}; \ 
{s\epsilon\over q}, \ 
{\alpha\over u},
{\beta \over u},
{\gamma \over u},
{\delta \over u}\big)\over
\big( {s\epsilon\over u^2}\big)_\infty
\big( {\alpha\epsilon\over u},
{\beta \epsilon\over u},
{\gamma\epsilon\over u},
{\delta \epsilon\over u}\big)_\infty
W \big( {\epsilon\over  u^2};\ \epsilon, 
{q\over \alpha u},
{q\over \beta u},
{q\over \gamma u},
{q\over \delta u}\big)} \Bigg] . \cr}
$$
Equating (4.6) and (4.7), and writing
$$
\leqalignno{
& G (\alpha, \beta, \gamma , \delta, \epsilon, u) := {1\over u}\ 
{\big( {\epsilon q\over u^2}\big)_\infty
\big( 
{s\epsilon\over \alpha u},
{s\epsilon\over \beta  u},
{s\epsilon\over \gamma  u},
{s\epsilon\over \delta u}\big)_\infty\over
\big(s\epsilon/ u^2\big)_\infty
\big(
{\alpha\epsilon\over u},
{\beta\epsilon\over u},
{\gamma \epsilon\over u},
{\delta \epsilon\over u}\big)_\infty} & (4.8)\cr
\times & 
W \big( {s\epsilon\over u^2 q} ; \ 
{s\epsilon\over q},
{\alpha\over u},
{\beta \over u},
{\gamma \over u},
{\delta \over u}\big)
W \big( \epsilon u^2;\ \epsilon, 
{qu\over \alpha},
{qu\over \beta},
{qu\over \gamma},
{qu\over \delta}\big)\ , \cr}
$$
we obtain the identity
$$
\leqalignno{
& G (\alpha, \beta, \gamma, \delta , \epsilon, u) - 
 G \big(\alpha, \beta, \gamma, \delta , \epsilon, 1/ u\big ) 
 \ = \ \big(1/ u - u\big) \ 
{(1-s\epsilon^2/ q^2 )\over
(1-s\epsilon/ q^2)}
 \cr
\times \ & 
{\big
(\alpha \beta \epsilon,
\alpha \gamma \epsilon,
\alpha \delta \epsilon,
\beta\gamma \epsilon,
\beta\delta \epsilon,
\gamma\delta \epsilon,
\epsilon q/ u^2,
\epsilon q u^2\big)_\infty\over
\big(\alpha\epsilon/ u ,
\alpha\epsilon u,
\beta \epsilon/ u,
\beta \epsilon u,
\gamma\epsilon/ u,
\gamma\epsilon u,
\delta\epsilon/ u,
\delta\epsilon u\big)_\infty }\cr}
$$
which is a  $ q $-analogue of Masson's generalization of Dougall's
theorem [18].  $ \epsilon = 1 $  gives a  $ q $-analogue of 
Dougall's theorem (see [3]).

We can recover from (4.9) the following identity which we had 
obtained in our earlier paper ([7], (31), p.723) for
$ \epsilon = 1 $,  $ s = q^m $;  $ m = 1,2,\ldots , $ i.e.
$$
\leqalignno{
& \big({q\over \alpha}\big)^{m-3}
\big( {1\over u} - u \big)
\Big[ u^{m-2} \Pi_1 (u) \Pi_2 \big({1\over u}\big) -
u^{2-m} \Pi_1 \big({1\over u}\big) \Pi_2 (u)\Big] \cr
& = ( \alpha \beta q^{-1}, 
\alpha\gamma q^{-1},
\alpha \delta q^{-1})_\infty
\big( {q^2\over \alpha\beta},
 {q^2\over \alpha\gamma },
 {q^2\over \alpha\delta }\big)_\infty
 \big(u^2, {1\over u^2}\big)_\infty \cr}
$$
for  $ m = 1,2, \ldots$ and $ |u| = 1 $,  where
$$
\leqalignno{
\Pi_1(u) \ = \ & \big( 
{q\over \alpha u},
{q\over \beta  u},
{q\over \gamma  u},
{q\over \delta  u}\big)_\infty \cr
\Pi_2(u) \ = \ & \big( 
{\alpha\over u},
{\beta\over  u},
{\gamma \over u},
{\delta\over  u}\big)_\infty\ . \cr}
$$
\vfill\eject

\centerline {\bf References }
\baselineskip=12pt
\frenchspacing
\vglue .25in
\item{1. }
R. Askey and J. Wilson,
Some basic hypergeometric orthogonal polynomials that generalize
Jacobi polynomials,
{\it Memoirs Amer. Math. Soc.} {\bf 319} (1985) 1--55.
\medskip

\item{2. }
R. Askey and J. Wimp,
Associated Laguerre and Hermite polynomials,
{\it Proc. Roy. Soc. Edinburgh} Sect. A 96 (1984), 15--37.
\medskip

\item{3. }
W.N. Bailey,
{\it Generalized Hypergeometric Series,}
Cambridge Univ. Press, London, 1935.
\medskip

\item{4. }
J. Bustoz and M.E.H. Ismail,
The associated ultraspherical polynomials and their $ q $-analogues,
{\it Canad. J. Math.} {\bf 34} (1982), 718--736.
\medskip

\item{5. } G. Gasper and M. Rahman,
{\it Basic Hypergeometric Series,}
Cambridge Univ. Press, Cambridge, 1990.
\medskip

\item{6. }
W. Gautschi,
Computational aspects of three-term recurrence relations,
{\it SIAM Rev.} {\bf 9} (1967), 24--82.
\medskip

\item{7. } D.P. Gupta and D.R. Masson,
Exceptional  $ q $-Askey-Wilson polynomials and continued fractions,
{\it Proc. A.M.S.} {\bf 112} (1991), 717--727.
\medskip

\item{8. } D.P. Gupta, M.E.H. Ismail and D.R. Masson,
Associated continuous Hahn polynomials,
{\it Canad. J. of Math.} {\bf 43} (1991), 1263--1280.
\medskip

\item{9. } D.P. Gupta, M.E.H. Ismail and D.R. Masson,
Contiguous relations, Basic Hypergeometric functions and
orthogonal polynomials II, Associated big $ q $-Jacobi polynomials,
{\it J. of Math. Analysis and Applications}  {\bf 171} (1992), 477--497.
\medskip

\item{10. } 
D.P. Gupta and D.R. Masson,
Watson's basic analogue of Ramanujan's Entry 40 and its
generalization,
{\it SIAM J. Math. Anal.}, to appear.
\medskip

\item{11. } M.E.H. Ismail, J. Letessier, and G. Valent,
Linear birth and death models and associated Laguerre polynomials,
{\it J. Approx. Theory} {\bf 56} (1988), 337--348.
\medskip

\item{12. } M.E.H. Ismail, J. Letessier, and G. Valent,
Quadratic birth and death processes and associated continuous
dual Hahn polynomials,
{\it SIAM J. Math. Anal.} {\bf 20} (1989), 727--737.
\medskip

\item{13. } M.E.H. Ismail, J. Letessier, G. Valent and J. Wimp,
Two families of associated Wilson polynomials,
{\it Can. J. Math.} {\bf 42} (1990), 659--695.
\medskip

\item{14. } M.E.H. Ismail and C.A. Libis,
Contiguous relations, basic hypergeometric functions and orthogonal
polynomials I,
{\it J. Math. Anal. Appl.} {\bf 141} (1989), 349--372.
\medskip

\item{15. } M.E.H. Ismail and D.R. Masson,
Two families of orthogonal polynomials related to Jacobi
polynomials,
{\it Rocky Mountain J. Math.} {\bf 21} (1991), 359--375. 
\medskip

\item{16. } M.E.H. Ismail and M. Rahman,
Associated Askey-Wilson polynomials,
{\it Trans. Amer. Math. Soc.} {\bf 328} (1991), 201--239.
\medskip

\item{17. }
W.B. Jones and W.J. Thron, {\it Continued Fractions: Analytic
Theory and Applications,} 
Addison-Wesley, Reading, Mass., 1980.
\medskip

\item{18. } D.R. Masson, 
Associated Wilson polynomials, 
{\it Constructive Approximation} {\bf 7} (1991), 521--534.
\medskip

\item{19. }
\bibline ,
Wilson polynomials and some continued fractions
of Ramanujan,
{\it Rocky Mountain J. of Math.} {\bf 21} (1991), 489--499.
\medskip

\item{20. }
\bibline ,
The rotating harmonic oscillator eigenvalue problems, I.
Continued fractions and analytic continuation,
{\it J. Math. Phys.} {\bf 24} (1983), 2074--2088.
\medskip

\item{21. }
\bibline,
Difference equations, continued fractions, Jacobi Matrices and
orthogonal polynomials, In:
{\it Non-linear numerical methods and Rational Approximation}
(A. Cuyt, ed.) Dordrecht, Reidel, 1988, 239--257.
\medskip

\item{22. }
J.A. Wilson,
Hypergeometric series, recurrence relations and some new orthogonal
polynomials,
Ph.D. diss., University of Wisconsin, Madison, 1978.
\medskip

\item{23. }
J. Wimp,
Explicit formulas for the associated Jacobi polynomials and
some applications,
{\it Canad. J. Math.} {\bf 39} (1987), 983--1000.
\end